\documentclass[11pt, oneside]{amsart}   	
\usepackage{geometry}                		
\usepackage{graphicx}				
\usepackage{amssymb}


\newtheorem{proposition}{Proposition}[section] 
\newtheorem{theorem}{Theorem}[section]
\newtheorem{definition}{Definition}[section]
\newtheorem{corollary}{Corollary}[section]
\newtheorem{lemma}{Lemma}[section]

\newtheorem*{theorem*}{Theorem}
\newtheorem{remark}{Remark}[section]
\newtheorem{example}{Example}[section]

\newcommand{\hf}{{\mathcal H}}

\newcommand{\fb}{{\mathcal B}}

\newcommand{\D}{{\mathcal D}}

\newcommand{\leb}{{\mathcal L}}
\newcommand{\M}{{\mathcal{M}}}

\newcommand{\Q}{{\mathcal Q}}
 \newcommand{\R}{{\mathcal R}}

\newcommand{\rn}{{\mathbb{R}}}

\newcommand{\z}{{\mathbb{Z}}}

\newcommand{\essint}{{\rm{int}}}
\newcommand{\esscl}{{\rm{cl}}}

\newcommand{\equalmm}{\;{= \hskip-3.4mm^{^\mu}}\;}

\newcommand{\re} {\hskip 1 pt{\hbox to 10.8 pt{\hfill \vrule
height7pt width0.4pt depth0pt \hbox{\vrule height0.4pt
width7.6pt depth0pt} \hfill}}}

\newcommand{\spt}{\operatorname{spt}}

\newcommand{\ree}{\mbox{\Large$\llcorner$}}

\newcommand{\avint}{{\mathchoice {\kern1ex\vcenter{\hrule height.4pt
width 6pt
depth0pt} \kern-9.7pt} {\kern1ex\vcenter{\hrule height.4pt width 4.3pt
depth0pt}
\kern-7pt} {} {} }\int}

\numberwithin{equation}{section}    

\makeatletter
\@namedef{subjclassname@2020}{%
  \textup{2020} Mathematics Subject Classification}
\makeatother

\title[Superdensity w.r.t. a Radon measure on $\rn^n$]{Superdensity with respect to \\ a Radon measure on $\rn^n$}  
\author{Silvano Delladio}

\begin{document}


\begin{abstract}
We introduce and investigate superdensity and the density degree of sets with respect to a Radon measure on $\rn^n$. Some applications are provided. In particular we prove a result on the approximability of a set by closed subsets of small density degree and a generalization of Schwarz's theorem on cross derivatives.
   \end{abstract}

  \maketitle

AFFILIATION: 
{\it 
University of Trento, Department of Mathematics via Sommarive 14, 38123 Trento, Italy}

\vskip3mm
{E-MAIL: {\it silvano.delladio@unitn.it}}

\vskip3mm
{ORCID: 0000-0002-5028-9573}

\vskip7mm

\section{Introduction}

Let us consider a Radon outer measure $\mu$ on $\rn^n$ and a $\mu$ measurable set $ E\subset\rn^n$. Then a celebrated result (cf. \cite[Cor.2.14]{Mattila1995}) states that for $\mu$ almost all $x\in E$ the set $E$ is $\mu$-dense at $x$, i.e.,
\begin{equation}
\label{eq:mdensitypointcond}
\lim_{r\to 0+} \frac{\mu (B_r(x)\cap E)}{\mu (B_r(x))}=1, \text{ which is equivalent to }
\lim_{r\to 0+} \frac{\mu (B_r(x)\setminus E)}{\mu (B_r(x))}=0,
\end{equation}
where $B_r(x)$ denotes the open ball in $\rn^n$, with centre $x$ and radius $r$. If the condition \eqref{eq:mdensitypointcond} is verified, then we can pose the problem of defining a number $d_E^\mu (x)$ that exactly quantifies the density of $E$ (w.r.t. $\mu$) at $x$. A natural way (not the only way, certainly!) to solve this problem is as follows:
\begin{itemize} 
\item First we say that $x$ is an $h$-superdensity point of $E$ (w.r.t. $\mu$) if $h\in [0,+\infty)$ and $\frac{\mu (B_r(x)\setminus E)}{\mu (B_r(x))}=o(r^h)$, as $r\to 0+$;
\item Then we define the density degree of $E$ (w.r.t. $\mu$) at $x$, denoted by $d_E^\mu (x)$, as the supremum of all $h\in [0,+\infty)$ such that $x$ is an $h$-superdensity point of $E$.
\end{itemize}
 In our previous work we have obtained a number of results concerning superdensity with respect to the Lebesgue outer measure $\leb^n$ and the purpose of the present paper is to generalize some of these results. 
 
 \medskip
 In this introduction we want to summarize the most significant parts of the paper. 
 Section \ref{sect:superdense} is devoted to prove some properties of the operator $b^{\mu,h}:2^{\rn^n}\to 2^{\rn^n}$ (with $h\in [0,+\infty)$) defined as follows
 $$
 b^{\mu,h} (A):=
 \left\{
 x\in\spt\mu \,\bigg\vert\,
 \limsup_{r\to 0+}\frac{\mu(B_r(x)\cap A)}{\mu(B_r(x)) r^h} >0
 \right\}\qquad (A\subset\rn^n).
 $$
 Roughly speaking, $b^{\mu,h} (A)$ is the set of all $x\in\spt\mu$ such that the relative size of $A$ in $B_r(x)$ is asymptotically larger than  $r^h$ (as $r\to 0+$).
 In Proposition \ref{prop:mainPropOfbmuh} we find that  $b^{\mu,h}$ is a base operator, i.e.,  $b^{\mu,h} (\emptyset)=\emptyset$ and 
 $$
 b^{\mu,h} (A\cup B)= b^{\mu,h}(A)\cup b^{\mu,h} (B) 
 $$
 for all $A,B\in 2^{\rn^n}$. Moreover, if $A^{\mu , h}$ denotes the set of all $h$-superdensity points of $A$ (w.r.t. $\mu$), then
 $$
 A^{\mu , h}\cup (\spt\mu)^c=[b^{\mu,h}(A^c)]^c.
 $$
  Hence $b^{\mu,h}$ determines a topology $\tau_{b^{\mu,h}}$ on $\rn^n$ which is finer than the ordinary Euclidean topology and such that 
  $$
  A\in \tau_{b^{\mu, h}}\text{ if and only if }A\cap \spt\mu \subset A^{\mu , h}, 
  $$
  cf. Proposition \ref{prop:ProprBaseTop}. 
  There are two main results in this paper. The first one, Theorem \ref{thm:GenProp3.2ofD32}, generalizes \cite[Prop.3.2]{D32}. It provides assumptions under which, in particular, the following property occurs (for any open set $\Omega\subset \rn^n$): For every $\varepsilon >0$ there exists an open set $A\subset\Omega$ such that $\mu (A) <\varepsilon$ and $A$ is so \lq\lq scattered\rq\rq\ that the inclusion $\Omega\cap\spt\mu\subset b^{\mu , h} (A)$ holds whenever $h$ exceeds a certain value which does not depend on $\varepsilon$. Here is the full statement:
  
  \medskip
  {\bf Theorem \ref{thm:GenProp3.2ofD32}.}\emph{ Let $\mu$ be non-trivial, i.e., $\spt\mu\not =\emptyset$. Suppose that there exist $C, p, q, \bar r \in (0, +\infty)$ such that $q\leq \min\{n , p\}$ and
\begin{equation*}
\frac{r^p}{C}
\leq  \mu (B_r(x))
\leq Cr^q
\end{equation*}
for all $x\in\spt \mu$ and $r\in (0,\bar r)$.  
 The following properties hold  for all $\varepsilon >0$ and $h>\frac{np}{q} - q$ (note that $\frac{np}{q} - q$ is non-negative): 
 \begin{itemize}
 \item[(1)] If $\Omega\subset\rn^n$ is a non-empty bounded open set, then there exists an open set $A\subset\Omega$ such that
 \begin{equation*}
 \mu (A)< \varepsilon,
 \qquad
 \Omega\cap \spt\mu \subset 
 b^{\mu, h} (A)
 \subset
\overline\Omega\cap \spt\mu.
 \end{equation*}
 In the special case when
 \begin{equation*}
 \partial\Omega\cap\spt\mu\subset b^{\mu,h}(\Omega),
 \end{equation*}
 the set $A$ can be chosen so that we have
 \begin{equation*}
  b^{\mu, h} (A)= \overline\Omega\cap \spt\mu.
 \end{equation*}
 \item[(2)] There is an open set $U\subset\rn^n$ satisfying
 $$
 \mu (U)<\varepsilon,\qquad
 b^{\mu,h}(U) 
 =\spt\mu.
 $$
 \end{itemize}
}

An example of application of  Theorem \ref{thm:GenProp3.2ofD32} to the Radon measure carried by a regular surface in $\rn^n$ is given in Subsection  \ref{subsect:2ndApplsOfFirstMainThm}. Another application is Proposition \ref{cor:ThinSetF}, which generalizes a property stated in \cite[Prop.5.4]{D34}. It provides a result on the approximability of a set by closed subsets of small density degree (w.r.t. $\mu$):

\medskip
{\bf Proposition \ref{cor:ThinSetF}.}\emph{
Let $\mu$ be non-trivial and assume that:
\begin{itemize}
\item[(i)] There exist $C, p, q, \bar r \in (0, +\infty)$ such that $q\leq \min\{n , p\}$ and
\begin{equation*}
\frac{r^p}{C}
\leq  \mu (B_r(x))
\leq Cr^q
\end{equation*}
for all $x\in\spt \mu$ and $r\in (0,\bar r)$;
\item[(ii)] It is given a non-empty bounded open set $\Omega\subset \rn^n$ with the following property: There exists an open bounded set $\Omega'\subset\rn^n$ such that $\Omega\subset\Omega'$ and $\partial\Omega' \cap \spt\mu\subset b^{\mu , h} (\Omega')$ for all $h>\overline m:= \frac{np}{q} - q$.
\end{itemize}
Then for all $H\in \left(0, \mu \left(\overline\Omega\right)\right)$ there exists a closed   subset $F$ of $\overline\Omega$ such that $\mu (F) > H$ and $d_F^\mu (x)\leq\overline m$ at $\mu$-a.e. $x$.
}

\medskip
The second main result generalizes the classical Schwarz theorem on cross derivatives (cf. Remark \ref{rem:ThmReducesToSchwarz} below). Here is the statement:

\medskip
{\bf Theorem \ref{thm:gengenschwarz}.}\emph{
  Let us consider an open set $\Omega\subset\rn^n$, $f, G, H\in C^1(\Omega)$, a couple of integers $p,q$ such that $1\leq p< q\leq n$ and $x\in\rn^n$. Assume that:
  \begin{itemize}
  \item[(i)] For $i=p, q$, the $i$-th distributional derivative of $\mu$ is a Borel real measure on $\rn^n$ also denoted $D_i\mu$, so that we have $D_i\mu (\varphi)=-\int D_i\varphi\, d\mu = \int\varphi\, d(D_i\mu)$, for all $\varphi\in C_c^1(\rn^n)$;
  \vskip2mm
  \item[(ii)] $x\in\Omega\cap A^{\mu , 1}$, where $A:=\{ y\in\Omega\,\vert\, (D_pf (y), D_qf(y)) =(G(y),H(y))\}$ (in particular $x\in\spt\mu$);
   \vskip2mm
  \item[(iii)] $\lim_{\rho\to 1-}\sigma (\rho)=1$, where $\sigma (\rho):= \liminf_{r\to 0+}\frac{\mu (B_r(x))}{\mu(B_{\rho r}(x))}$ (note that $\sigma$ is decreasing);
   \vskip2mm
  \item[(iv)] For $i=p,q$, one has $\lim_{r\to 0+}\frac{\vert D_i\mu\vert (B_r(x))}{r\mu (B_r(x))}=0$ (where $\vert D_i\mu\vert$ denotes the total variation of $D_i\mu$). 
  \end{itemize}
  Then $D_pH (x) = D_qG(x)$.
}

\medskip
Among the results obtained in our previous work are several of the same kind as Theorem \ref{thm:gengenschwarz}, in the special case $\mu=\leb^n$. They were then applied to describe the fine properties of sets of solutions of differential identities under assumptions of non-integrability. The simplest example that we can mention is $Df=F$, with $f\in C^1(\rn^2)$ and $F=(F_1,F_2)\in C^1(\rn^2 , \rn^2)$ such that $D_1F_2 (x)\not = D_2 F_1(x)$ for every $x\in\rn^3$. If we recall that $D_1\leb^2= D_2\leb^2=0$ and apply Theorem \ref{thm:gengenschwarz} with
$$
n=2,\quad
\Omega=\rn^2, \quad
G=F_1,\quad
H=F_2,\quad
p=1,\quad
q=2, \quad
\mu=\leb^2,
$$
then we conclude that $A^{\leb^2 , 1}=\emptyset$, regardless of $f$,  even though there are functions $f$ such that $\leb^2 (A)>0$ (cf. \cite[Theorem 2.1]{D26}). In particular, the density degree of $A$ (w.r.t. $\leb^2$) is less than or equal to $1$ everywhere and  this gives us fairly accurate information about the fine structure of $A$. Similar arguments have been used, for example:
\begin{itemize} 
\item In \cite{D37}, to prove that, given a $C^1$ smooth $n$-dimensional submanifold $M$ of $\rn^{n+m}$ and a non-involutive $C^1$ distribution $\D$ of rank $n$ on $\rn^{n+m}$, the tangency set of $M$ with respect to $\D$ can never be too dense. 
\item In \cite{D44, D45}, to obtain results about low density of the set of solutions of the differential identity $G(D)f=F$, for certain classes of linear partial differential operators $G(D)$, under assumptions of non-integrability on $F$.
\end{itemize}
In connection with the results in \cite{D26} and \cite{D37}, we would like to mention the paper \cite{AMS2022} on the structure of tangent currents to smooth distributions. The application of superdensity used in \cite{D50}  is a first successful attempt to extend the theory developed so far for the Lebesgue measure to other contexts (tangency of generalized surfaces as considered in \cite{AMS2022}). At the same time, it gives us reason to believe that it is interesting to continue working on generalisation. It is in this sense that the present work, which provides a superdensity theory for Radon measures on $\rn^n$, should be understood.


\section{Basic notation and notions}

\subsection{Basic notation}
The Lebesgue outer measure on $\rn^n$ and the  $s$-dimensional Hausdorff outer measure on $\rn^n$ are denoted by $\leb^n$ and $\hf^s$, respectively. The $i$-th partial derivative, either classical or distributional, will be denoted by $D_i$. The ordinary topology of $\rn^n$ is denoted by $\tau (\rn^n)$.
The $\sigma$-algebra generated by $\tau (\rn^n)$ is denoted by $\fb (\rn^n)$. A member of $\fb (\rn^n)$ is called \emph{Borel set}.
$B_r(x)$ is the open ball in $\rn^k$, with centre $x$ and radius $r$ ($k$ does not appear in the notation as its value will be made clear from the context). 
The family of all Radon outer measures on $\rn^n$ is denoted by $\R$. If $\mu\in\R$ then $\M_\mu$ is the $\sigma$-algebra of all $\mu$ measurable sets. When two subsets $A$ and $B$ of $\rn^n$ are equivalent with respect to $\mu\in\R$, i.e., $\mu (A\setminus B)= \mu (B\setminus A)=0$, we write $A\equalmm B$. Observe that if $A\equalmm B$ and $B\in\M_\mu$, then $A\in\M_\mu$.
If $\mu\in\R$ then $\spt \mu$ denotes the support of $\mu$, that is the smallest closed set $F\subset\rn^n$ such that $\mu (\rn^n\setminus F)=0$. Hence
\begin{equation}
\label{eq:mu(sptmuc)=0}
\mu ( \rn^n\setminus\spt\mu)=0
\end{equation} 
and 
\begin{equation}
\label{eq:ComplSpt=Nmu}
 \rn^n\setminus\spt\mu =
 \{x\in\rn^n\,\vert\,
\mu(B_r(x))=0 \text{ for some } r>0
\},
\end{equation}
cf. \cite[Def.1.12]{Mattila1995}. The total variation of a Borel real measure $\lambda$ on $\rn^n$ is denoted by $\vert\lambda\vert$ (cf. \cite[Def.1.4]{AFP2006}).
 


\subsection{superdensity}
\label{subsect:Superdensity} The following definition has been introduced in \cite{D50} and generalizes the notion of $m$-density point (cf. \cite{D26, D27, D28}).
 
 \begin{definition}
 \label{def:SuperdensWRTmu}
 Let $\mu\in\R$, $h\in [0,+\infty)$ and $E\subset\rn^n$. Then $x\in\rn^n$ is said to be an $h$-superdensity point of $E$ w.r.t. $\mu$ if $x\in\spt\mu$ and $\mu(B_r(x)\setminus E)=\mu (B_r(x))\, o(r^h)$, as $r\to 0+$. The set of all $h$-superdensity points of $E$ w.r.t. $\mu$ is denoted by $E^{\mu,h}$.
 \end{definition}
 
 \begin{remark}
 \label{rem:PropSuperd}
 Let $\mu\in\R$, $h\in [0,+\infty)$ and $E, F\subset\rn^n$. Then it can easily be verified that the following properties hold true:
 \begin{itemize}
  \item[(1)] If $\mu=\leb^n$ then the set of all $h$-superdensity points of $E$ w.r.t. $\mu$ coincides with the set of all $(n+h)$-density points of $E$, i.e., $E^{\leb^n, h}=E^{(n+h)}$.
 \item[(2)] $E^{\mu,h_2}\subset E^{\mu,h_1}$, whenever $0\leq h_1\leq h_2<+\infty$.
  \item[(3)] $(E\cap F)^{\mu , h} = E^{\mu , h}\cap F^{\mu , h}$.
  \item[(4)] If $E,F\in\M_\mu$ and $E\equalmm F$, then $E^{\mu , h} = F^{\mu , h}$.  In particular, this equality occurs whenever $E\in\M_\mu$ has finite measure and $F$ is a \lq\lq Borel envelope\rq\rq\ of $E$ (that is $F\in\fb (\rn^n)$, $F\supset E$ and $\mu (F)=\mu (E)$).
   \item[(5)] If $E\in\M_\mu$, then $E^{\mu , 0}\equalmm E$ (cf. \cite[Cor.2.14]{Mattila1995}).
   \item[(6)] Let $E$ be open. Then $E\subset E^{\mu , h}$ and the inclusion can be strict, e.g., for $\mu =\leb^n$ and $E=B_r(x)\setminus\{x\}$ one has $E^{\mu , h}=B_r(x)$.  
   \item[(7)] $E^{\mu , h}\subset\spt\mu\cap \overline E$. In particular, if $E$ is closed then $E^{\mu , h}\subset \spt\mu\cap E$.
   \item[(8)] $E^{\mu\ree E , h}=\spt\mu$.
   \end{itemize}
 \end{remark}
 
 \begin{remark}
 Recall from \cite[Lemma 4.1]{D26} that if $E$ is a locally finite perimeter subset of $\rn^n$ (cf. \cite[Sect.3.3]{AFP2006}) then $\leb^n(E\setminus E^{\leb^n , \frac{n}{n-1}})=0$.
 \end{remark}
 
 \subsection{Base operators}
 Let us recall from \cite[Ch.1]{LMZ1986} that a \emph{base operator} on a set $X$ is a map $b:2^X\to 2^X$ such that $b (\emptyset)=\emptyset$ and $b (A\cup B) = b (A)\cup b(B)$ for all $A,B\in 2^X$. Any base operator $b$ is obviously monotone and determines a topology on $X$ that is defined as follows
 $$
 \tau_b:=\left\{
 A\in 2^X\big\vert
 \, b(X\setminus A)\subset X\setminus A
 \right\}.
 $$ 
 It turns out that $\tau_b$ is the finest topology $\tau$ on $X$ such that, for all $A\subset X$, the closure of $A$ w.r.t. $\tau$ contains $b(A)$. If $X=\rn^n$ and $b(A)$ denotes the ordinary closure of $A\subset\rn^n$, then $b$ is a base operator and $\tau_b=\tau (\rn^n)$.

\section{Superdensity w.r.t. the measure carried by a regular surface}
\label{sect:ExampleC1Map}

Let $G$ be a bounded open subset of $\rn^k$ and consider $\varphi \in C^1 (\rn^k, \rn^{n})$ such that $\varphi\vert_G$ is an imbedding ($k\leq n$). In particular
\begin{equation}
\label{eq:detDphi2>0}
J\varphi (y):=
\left[\det [(D\varphi)^t\times (D\varphi)] (y)\right]^{1/2} >0
\end{equation}
for all $y\in G$.
We observe that $\hf^k\ree \varphi (G)\in\R$.

\medskip
We will prove:

\begin{proposition}
\label{prop:phi(E(k+h))=E(Hkphih)}
If $E\subset\rn^{n}$ and $h\in [0,+\infty)$, then
\begin{equation*}
E^{\hf^k\ree \varphi (G) , h}\cap \varphi (G)
=\varphi\left(
[\varphi^{-1}(E)]^{(k+h)}\cap G
\right).
\end{equation*}
\end{proposition}

\begin{remark}
\label{rem:P-1Ek+h=PG-1Ek+h}
From (3) and (6) in Remark \ref{rem:PropSuperd}, it follows that
\begin{equation*}
[\varphi^{-1}(E)]^{(k+h)}\cap G
=
[\varphi^{-1}(E)\cap G]^{(k+h)}\cap G
=
[(\varphi\vert_G )^{-1}(E)]^{(k+h)}\cap G.
\end{equation*}
\end{remark}

In the proof of Proposition \ref{prop:phi(E(k+h))=E(Hkphih)}, we will need the following easy corollary of  \cite[Ch.VIII, Th.3.3]{Lang1987}.

\begin{lemma}
\label{lem:minLuu>0}
Let $L$ be a real symmetric matrix of order $k$ such that $\det L\not =0$ and $(Lv)\cdot v\geq 0$ for all $v\in\rn^k$. Then 
$\min \{(Lu)\cdot u\,\vert  \, u\in\rn^k,\, \vert u\vert=1\} >0$.
\end{lemma}

\medskip
\emph{Proof of Proposition \ref{prop:phi(E(k+h))=E(Hkphih)}.}
Let us consider an arbitrary $y\in G$. We have to prove that
\begin{equation*}
\varphi (y)\in E^{\hf^k\ree \varphi (G) , h}
\text{ if and only if }
y\in [\varphi^{-1}(E)]^{(k+h)}
\end{equation*}
namely, setting for simplicity $\mu:=\hf^k\ree \varphi (G)$,
\begin{equation}
\label{eq:goal}
\lim_{r\to 0+}
\frac{\mu (B_r(\varphi(y))\cap E^c)}{\mu (B_r(\varphi(y)))r^h}=0
\text{ if and only if }
\lim_{r\to 0+}
\frac{\leb^k (B_r(y)\cap [\varphi^{-1}(E)]^c)}{r^{k+h}}=0.
\end{equation}
To this end we observe that
\begin{equation}
\label{eq:deltaphi1}
\begin{split}
\varphi (z)-\varphi (y)
&= \int_0^1(D\varphi)(y+t(z-y)) (z-y) \, dt
\end{split}
\end{equation}
for all $z\in\rn^k$. If $\Vert \cdot \Vert$ denotes the Hilbert-Schmidt norm of matrices and we define
\begin{equation}
\label{def:Kandm1}
K:=\{
z\in\rn^k\,\vert\, \text{dist} (z, G)\leq 1
\}, \quad
m_1 :=\max_{z\in K}\Vert (D\varphi) (z)\Vert > 0
\end{equation}
then \eqref{eq:deltaphi1} yields
\begin{equation}
\label{eq:deltaphileqMr}
\begin{split}
\vert\varphi (z)-\varphi (y)\vert
&\leq r\int_0^1\Vert (D\varphi)(y+t(z-y))\Vert \, dt
\leq m_1 r
\end{split}
\end{equation}
for all $z\in\overline{B_r(y)}$ with $r\in (0,1]$.
 Furthermore, by \eqref{eq:deltaphi1}, we have
 \begin{equation*}
\begin{split}
\varphi (z)-\varphi (y)
&=(D\varphi)(y)  (z-y) +
 \int_0^1 [(D\varphi)(y+t(z-y))-(D\varphi)(y)]  (z-y) \, dt
\end{split}
\end{equation*}
for all $z\in\rn^k$. Hence, for all $r>0$ and $z\in\partial B_r(y)$, we obtain
\begin{equation}
\label{eq:deltaphigeqmr-sr}
\begin{split}
\vert\varphi (z)-\varphi (y)\vert
&\geq
 \left[\left([(D\varphi)^t\times (D\varphi)] (y) (z-y)\right)\cdot (z-y)\right]^{1/2}\\
 & -
r \int_0^1 \Vert (D\varphi)(y+t(z-y))-(D\varphi)(y)\Vert \, dt\\
&\geq 2 m_0 r - \sigma_r  r
\end{split}
\end{equation}
where
\begin{equation*}
m_0:= \frac{1}{2}\left[\min \left\{\left([(D\varphi)^t\times (D\varphi)] (y)u\right)\cdot u\,\vert  \, u\in\rn^k,\, \vert u\vert=1\right\}\right]^{1/2}
\end{equation*}
and
\begin{equation*}
\sigma_r:= \max_{z\in\overline{B_{r}(y)}}\Vert (D\varphi) (z) - (D\varphi) (y)\Vert.
\end{equation*}
Observe that:
\begin{itemize}
\item Since $\varphi$ is of class $C^1$, then 
\begin{equation}
\label{eq:sigmainfinitesimal}
\lim_{r\to 0+}\sigma_r = 0;
\end{equation}
\item \eqref{eq:detDphi2>0} and Lemma \ref{lem:minLuu>0} with $L=[(D\varphi)^t\times (D\varphi)] (y)$, yield 
\begin{equation}
\label{eq:m>0}
m_0>0.
\end{equation}
\end{itemize}
From \eqref{eq:deltaphigeqmr-sr}, \eqref{eq:sigmainfinitesimal} and \eqref{eq:m>0}, it follows that
\begin{equation}
\label{eq:deltaphigeqmr}
\vert\varphi (z)-\varphi(y)\vert \geq m_0 r,
\text{ for all $z\in\partial B_r (y)$,}
\end{equation}
provided $r$ is small enough.
Now, by \eqref{eq:deltaphileqMr} and \eqref{eq:deltaphigeqmr}, we obtain
\begin{equation*}
\varphi (G)\cap B_{m_0 r}(\varphi (y))
\subset \varphi (B_r(y))
\subset \varphi (G)\cap B_{m_1 r}(\varphi (y)),
\end{equation*}
provided $r$ is small enough. Recalling also the area formula (cf. \cite[Cor. 5.1.13]{KrantzParks2008}), it follows that this set of inequalities holds for $r$ small enough:

\begin{equation*}
\mu \left(\varphi(B_{r/m_1 }(y))\right)
\leq \mu \left(B_{r}(\varphi(y))\right)
\leq \mu \left(\varphi(B_{r/m_0}(y))\right)
\end{equation*}
\begin{equation*}
\mu \left(\varphi(B_{r/m_1 }(y))\cap E^c\right)
\leq \mu \left(B_{r}(\varphi(y))\cap E^c\right)
\leq \mu \left(\varphi(B_{r/m_0}(y))\cap E^c\right)
\end{equation*}
\begin{equation*}
\frac{J\varphi(y)}{2} \leb^k (B_{r }(y))\leq
\mu \left(\varphi(B_{r }(y))\right)
=\int_{B_{r }(y)}J\varphi\, d\leb^k\leq 2J\varphi(y)\leb^k (B_{r }(y))
\end{equation*}
\begin{equation*}
\begin{split}
\frac{J\varphi(y)}{2} \leb^k (B_{r}(y)\cap [\varphi^{-1}(E)]^c)
&\leq\mu \left(\varphi(B_{r}(y))\cap E^c\right)
=\int_{B_{r}(y)\cap \varphi^{-1}(E)^c}J\varphi\, d\leb^k\\
&\leq 2 J\varphi(y) \leb^k (B_{r}(y)\cap [\varphi^{-1}(E)]^c).
\end{split}
\end{equation*}
Hence, the statement \eqref{eq:goal} follows easily.
\qed

 \section{Base operators associated to a Radon measure}
\label{sect:superdense}


 \begin{proposition}
 \label{prop:mainPropOfbmuh}
 Let $\mu\in\R$, $h\in [0, +\infty)$ and consider the operator $b^{\mu,h}:2^{\rn^n}\to 2^{\rn^n}$ defined as follows (recall \eqref{eq:ComplSpt=Nmu})
 $$
 b^{\mu,h} (A):=
 \left\{
 x\in\spt\mu \,\bigg\vert\,
 \limsup_{r\to 0+}\frac{\mu(B_r(x)\cap A)}{\mu(B_r(x)) r^h} >0
 \right\}\qquad (A\subset\rn^n).
 $$
 Then
 \begin{itemize}
 \item[(1)] $b^{\mu,h}(A)\subset \spt\mu$, for all $A\in 2^{\rn^n}$;
 \item[(2)] $A^{\mu , h}\cup (\spt\mu)^c=[b^{\mu,h}(A^c)]^c$, for all $A\in 2^{\rn^n}$;
 \item[(3)] $A\cap\spt\mu \subset b^{\mu , h} (A)$, for all $A\in\tau (\rn^n)$.
 \end{itemize}
 Moreover $b^{\mu,h}$ is a base operator, that is:
 \begin{itemize}
  \item[(4)] $b^{\mu,h} (\emptyset)=\emptyset$;
   \item[(5)] $b^{\mu,h} (A\cup B)= b^{\mu,h}(A)\cup b^{\mu,h} (B)$, for all $A,B\in 2^{\rn^n}$.
 \end{itemize}
 \end{proposition}

 \begin{proof}
 Statements (1), (2), (3) and (4) are trivial, while (5) follows easily by combining property (3) in Remark \ref{rem:PropSuperd} and (2).
 \end{proof}

 \begin{example}
 Given $\bar x\in\rn^n$, let $\delta_{\bar x}$ be the Dirac outer measure (on $\rn^n$) at $\bar x$. Then $\delta_{\bar x}\in\R$, $\spt\delta_{\bar x}= \{ \bar x\}$ and
 $$
 b^{\delta_{\bar x} , h} (A)=
 \begin{cases}
      \{  \bar x\} & \text{if $ \bar x\in A$} \\
      \emptyset & \text{if $ \bar x\not\in A$}
\end{cases}
 $$
 for all $h\in [0,+\infty)$ and $A\subset\rn^n$. Hence and recalling (2) of Proposition \ref{prop:mainPropOfbmuh} (or also simply by Definition \ref{def:SuperdensWRTmu}), we obtain
 $$
 A^{\delta_{\bar x} , h}= [  b^{\delta_{\bar x} , h} (A^c)]^c\cap \{\bar x\}
 =
 \begin{cases}
      \{\bar x\} & \text{if $ \bar x\in A$} \\
      \emptyset & \text{if $ \bar x\not\in A$}
\end{cases}
 $$
 for all $h\in [0,+\infty)$ and $A\subset\rn^n$. Moreover it is very easy to verify that $\tau_{b^{\delta_{\bar x} , h}} = 2^{\rn^n}$, for all $h\in [0,+\infty)$.
 \end{example}

The same arguments used in \cite[Prop.3.1]{D32} yield the following proposition.

\begin{proposition}
\label{prop:ProprBaseTop}
Let $\mu\in\R$ and $h\in [0, +\infty)$. The following facts hold:
\begin{itemize}
\item[(1)] $b^{\mu, h} (A)\in\M_\mu$, for all $A\in 2^{\rn^n}$. Hence $A^{\mu , h}\in \M_\mu$, for all $A\in2^{\rn^n}$.
\item[(2)] $A\in \tau_{b^{\mu, h}}$ if and only if $A\cap \spt\mu \subset A^{\mu , h}$. In particular $\tau (\rn^n)\subset \tau_{b^{\mu, h}}$.
\item[(3)] If $l\in [h,+\infty)$ then $b^{\mu, h} (A)\subset b^{\mu, l} (A)$, for all $A\in 2^{\rn^n}$. Hence $\tau_{b^{\mu, l}}\subset \tau_{b^{\mu, h}}$.
\end{itemize}
\end{proposition}

  The proof of Theorem \ref{thm:GenProp3.2ofD32} below is a non-trivial adaptation of the argument used to prove \cite[Prop.3.2]{D32}. We need to make a premise about lattices, which we include in the following remark. 
  
  \begin{remark}
  \label{rem:premretic}
  We consider three positive integers  $R, \beta, k$ and set $L_k:=(2R\beta^k)^n$. Let $P_1^{(k)},\ldots , P_{L_k}^{(k)}$
  be the points of the lattice $\Lambda_k:=(\beta^{-k}\z^n)\cap [-R,R)^n$ and define the corresponding cells (which we will simply call $k$-cells) as
  $$
  \Q_j^{(k)}
  :=P_j^{(k)} + [0,\beta^{-k})^n\qquad
  (j=1,\ldots,L_k).
  $$
  Observe that the $k$-cells form a partition of $[-R,R)^n$. Now let $S$ be an infinite subset of $[-R,R)^n$ and denote by $N_k$ the number of $k$-cells intersecting $S$. Obviously one has $N_k\leq N_{k+1}$ (for all $k\geq 1$) and $N_k\to +\infty$ (as $k\to +\infty$). Then we can easily find a countable family $\{ P_j\}\subset S$ such that the following property holds, for all $k\geq 1$: Each one of the $k$-cells intersecting $S$ contains one and only one point of $\{ P_1, P_2, \ldots, P_{N_{k}}\}$.
  
 Under the assumptions above, we finally define $\Lambda :=\cup_{k=1}^{+\infty}\Lambda_k$ and we say that $\{ P_j\}$ is a $\Lambda$-distribution of $S$.
  \end{remark}

\begin{theorem}
\label{thm:GenProp3.2ofD32}
Let $\mu\in\R$ be non-trivial, i.e., $\spt\mu\not =\emptyset$. Suppose that there exist $C, p, q, \bar r \in (0, +\infty)$ such that $q\leq \min\{n , p\}$ and
\begin{equation}
\label{eq:rpmuBrq}
\frac{r^p}{C}
\leq  \mu (B_r(x))
\leq Cr^q
\end{equation}
for all $x\in\spt \mu$ and $r\in (0,\bar r)$.  
 The following properties hold  for all $\varepsilon >0$ and $h>\frac{np}{q} - q$ (note that $\frac{np}{q} - q$ is non-negative): 
 \begin{itemize}
 \item[(1)] If $\Omega\subset\rn^n$ is a non-empty bounded open set, then there exists an open set $A\subset\Omega$ such that
 \begin{equation}
 \label{eq:(3.2)inD32}
 \mu (A)< \varepsilon,
 \qquad
 \Omega\cap \spt\mu \subset 
 b^{\mu, h} (A)
 \subset
\overline\Omega\cap \spt\mu.
 \end{equation}
 In the special case when
 \begin{equation}
 \label{eq:DOcapsptmsubb}
 \partial\Omega\cap\spt\mu\subset b^{\mu,h}(\Omega),
 \end{equation}
 the set $A$ can be chosen so that we have
 \begin{equation}
 \label{eq:(3.2)inD32bis}
  b^{\mu, h} (A)= \overline\Omega\cap \spt\mu.
 \end{equation}
 \item[(2)] There is an open set $U\subset\rn^n$ satisfying
 $$
 \mu (U)<\varepsilon,\qquad
 b^{\mu,h}(U) 
 =\spt\mu.
 $$
 \end{itemize}
\end{theorem}
 
 \begin{proof}
 First of all observe that, by  \eqref{eq:mu(sptmuc)=0} and \eqref{eq:rpmuBrq}, we have $\mu(\spt\mu)>0$ and 
 \begin{equation}
 \label{eq:mx=0}
 \mu(\{x\})=0, \text{ for all $x\in\spt\mu$}. 
 \end{equation}
 Hence $\spt\mu$ is a non-countable set. That said, we can proceed to prove (1) and (2).

 \medskip
 \underline{Proof of (1).} 
 If 
 \begin{equation}
 \label{eq:Ocapsptm=empty}
 \Omega\cap\spt\mu =\emptyset
 \end{equation} 
 holds, then:
 \begin{itemize}
\item The first statement is trivially verified with $A=\emptyset$. 
\item We have 
\begin{equation}
\label{eq:bO=emptyforallh}
b^{\mu,h}(\Omega)=\emptyset, \text{ for all $h\in (0,+\infty)$.}
\end{equation} 
For if this were not true, $x\in b^{\mu,h}(\Omega)$ would exist for a certain $h\in (0,+\infty)$ and this would imply $\mu((B_r(x)\setminus\{x\})\cap\Omega)>0$ for all $r>0$ (by \eqref{eq:mx=0}), which contradicts \eqref{eq:Ocapsptm=empty}. Now, in the special case when \eqref{eq:DOcapsptmsubb} holds, the equality \eqref{eq:bO=emptyforallh} yields $\partial\Omega\cap\spt\mu=\emptyset$ and it follows immediately from this that the second statement is also true.
\end{itemize}

Thus, we can assume that 
\begin{equation*}
\Omega\cap\spt\mu \not =\emptyset.
\end{equation*}

This assumption and \eqref{eq:ComplSpt=Nmu} (or \eqref{eq:rpmuBrq}) imply that there exists an open ball $B\subset\Omega$ such that $\mu(B)>0$, hence
\begin{equation*}
\mu (\Omega\cap\spt\mu)\geq
\mu (B\cap\spt\mu)=\mu (B)>0.
\end{equation*}

 From this fact and \eqref{eq:mx=0}, it follows that $\Omega\cap\spt\mu$ is a non-countable set.
   
 \medskip
 Now consider $\varepsilon >0$ and $h>\frac{np}{q} - q$. Define
 $$
 m:=\frac{(h+q)q}{p}.
 $$
 and observe that 
 \begin{equation}
 \label{eq:m>n}
 m>n, \text{ hence also }\frac{m}{q}>1. 
 \end{equation}
 Moreover let $R$ and $\beta$ be positive integers such that
 \begin{equation*}
 \Omega\subset [-R , R]^n
 \end{equation*}
 and
 \begin{equation}
 \label{eq:(3.4)inD32}
 \beta > \max
 \left\{
 (2^nR^n+1)^{\frac{1}{m-n}}\, ; \,
 \left(\frac{\varepsilon}{C}\right)^{1/q}+ n^{1/2}\, ;\,
 \left(\frac{\varepsilon}{C\bar r^q}\right)^{1/m}
 \right\}.
 \end{equation}
 For $k=1,2,\ldots$, we define 
 \begin{equation*}
 \rho_k:=
 \left(
 \frac{\varepsilon}{C\beta^{km}}
 \right)^{1/q}, \quad
 \Lambda_k :=(\beta^{-k}\z^n)\cap [-R,R)^n
 \end{equation*}
 and note that
 \begin{equation}
 \label{eq:rhokleqbr}
 \rho_k < \bar r
 \end{equation}
 by \eqref{eq:(3.4)inD32}. Then, by recalling Remark \ref{rem:premretic} and the notation therein, we can find a $\Lambda$-distribution $\{ P_j\}_{j=1}^\infty$ of $\spt\mu\cap \Omega$. We set (for $k=1,2,\ldots$)
 \begin{equation*}
 \label{eq:(3.5)inD32}
 \Gamma_k:=\{P_j\,\vert\, 
 1\leq j\leq N_k,\,
 B_{\rho_k}(P_j)\subset\Omega \},\quad
 A_k:=\bigcup_{P\in\Gamma_k} B_{\rho_k}(P), \quad
 A:=\bigcup_{k=1}^{+\infty}A_k
 \end{equation*}
 and observe that 
 \begin{equation}
 \label{eq:(3.6)inD32}
 \# (\Gamma_k)\leq N_k\leq
 L_k
 =2^nR^n\beta^{kn}.
 \end{equation}
 By  \eqref{eq:(3.4)inD32}, \eqref{eq:rhokleqbr}, \eqref{eq:(3.6)inD32} and assumption \eqref{eq:rpmuBrq}, we get
 \begin{equation*}
 \label{eq:(*)13_2809}
 \begin{split}
 \mu (A)
 &\leq \sum_{k=1}^{+\infty}\mu (A_k)
 \leq\sum_{k=1}^{+\infty}\sum_{P\in\Gamma_k}\mu(B_{\rho_k}(P))
 \leq C\, \sum_{k=1}^{+\infty}\#(\Gamma_k)\rho_k^q
 \leq \frac{2^nR^n\varepsilon}{\beta^{m-n}-1}
 < \varepsilon.
 \end{split}
 \end{equation*}
 Let us prove that
 \begin{equation}
 \label{eq:(3.7)inD32}
 \Omega\cap \spt\mu
 \subset b^{\mu, h}(A).
 \end{equation}
 To this end, consider $x\in\Omega\cap\spt\mu$ and chose $K_x>0$ such that
 \begin{equation*}
 \label{eq:(3.8)inD32}
 B_{\beta^{-K_x}}(x)\subset\Omega.
 \end{equation*}
 Obviously, for every $k\geq K_x +1$ there exists a $k$-cell containing $x$. This $k$-cell must also contain a point of $\{P_1, P_2, \ldots , P_{N_k}\}$, which we denote by $Q_k$ (cf. Remark \ref{rem:premretic}). Observe that
 \begin{equation*}
 \vert Q_k - x\vert
 \leq \beta^{-k}n^{1/2}.
 \end{equation*}
  Then, for all $k\geq K_x + 1$ and $y\in B_{\rho_{k}}(Q_k)$, we find (recalling \eqref{eq:(3.4)inD32} and \eqref{eq:m>n} too)
 \begin{equation*}
 \begin{split}
 \vert y-x\vert
 &\leq  \vert y-Q_k\vert +  \vert Q_k -x\vert
 <\rho_{k} + \beta^{-k}n^{1/2}
 = \left(
 \frac{\varepsilon}{C\beta^{km}}
 \right)^{1/q} + \beta^{-k}n^{1/2}\\
 &< \left[
 \left(
 \frac{\varepsilon}{C}
 \right)^{1/q}
 + n^{1/2}
 \right]\beta^{-k}
 <\beta^{-k+1}
 \leq\beta^{-K_x}.
 \end{split}
 \end{equation*}
 Thus
 \begin{equation}
 \label{eq:(**)2010p1}
 B_{\rho_{k}}(Q_k)
 \subset B_{\beta^{-k+1}}(x)\subset B_{\beta^{-K_x}}(x)\subset\Omega.
 \end{equation}
 In particular
 \begin{equation*}
 Q_k\in\Gamma_k
 \end{equation*}
 and hence
 \begin{equation}
 \label{eq:(*)2010p1}
 B_{\rho_{k}}(Q_k)
 \subset  A_{k}
 \subset A.
 \end{equation}
 From \eqref{eq:(**)2010p1} and \eqref{eq:(*)2010p1}, recalling \eqref{eq:rhokleqbr} and \eqref{eq:rpmuBrq} too, we obtain
 \begin{equation*}
 \begin{split}
 \mu\left(
 A  \cap B_{\beta^{-k+1}}(x)
 \right)
 \geq \mu\left(
 B_{\rho_{k}}(Q_k)
 \right)
 \geq \frac{\rho_{k}^p}{C}
 = \frac{\varepsilon^{p/q}\beta^{-kmp/q}}{C^{1+p/q}}.
 \end{split}
 \end{equation*}
 Hence, by \eqref{eq:rpmuBrq} and recalling the definition of $m$, we obtain (for $k$ large enough)
 \begin{equation*}
 \frac{ \mu\left(
 A  \cap B_{\beta^{-k+1}}(x)
 \right)}{ \mu\left(
  B_{\beta^{-k+1}}(x)
 \right) (\beta^{-k+1})^h}
 \geq
 \frac{\varepsilon^{p/q}\beta^{-kmp/q}}{C^{2+p/q}\beta^{(-k+1)q}\beta^{(-k+1)h}}
 =\frac{\varepsilon^{p/q}}{C^{2+p/q}\beta^{q+h}}
 \end{equation*}
 which shows that $x\in b^{\mu, h}(A )$ and concludes the proof of \eqref{eq:(3.7)inD32}. By recalling that
 \begin{itemize}
 \item $A\subset\Omega\subset\overline\Omega$,
 \item $b^{\mu , h}(\overline\Omega)\subset\spt\mu$ (cf.(1) in Proposition \ref{prop:mainPropOfbmuh}),
 \item $\overline\Omega$ is closed with respect to $\tau_{b^{\mu , h}}$ (cf.(2) in Proposition \ref{prop:ProprBaseTop}),
 \end{itemize}
 we can now complete the proof of \eqref{eq:(3.2)inD32}:
 \begin{equation}
 \label{eq:bAsubsClOintS}
 b^{\mu , h} (A)
 \subset b^{\mu , h}(\overline\Omega)
 =b^{\mu , h}(\overline\Omega)\cap\spt\mu
 \subset \overline\Omega\cap\spt\mu.
 \end{equation}
 Now assume that \eqref{eq:DOcapsptmsubb} holds. Then consider an open set $A'\subset\rn^n$ satisfying
 \begin{equation}
 \label{eq:A'contQmenoO}
 A'\supset [-R , R]^n\setminus\Omega, \quad
 \mu\left(
 A'\setminus ([-R , R]^n\setminus\Omega)
 \right) <\varepsilon-\mu (A).
 \end{equation}
 Observe that 
 \begin{equation}
 \label{eq:A'capOcontA'menoECC}
 A'\cap\Omega\subset
 A'\setminus ([-R , R]^n\setminus\Omega)
 \end{equation}
 and define
 \begin{equation}
 \label{eq:defA''}
 A''
 := A\cup (A'\cap \Omega),
 \end{equation}
 which is an open subset of $\Omega$. We shall prove that $A''$ satisfies \eqref{eq:(3.2)inD32} and \eqref{eq:(3.2)inD32bis}, that is
 \begin{equation}
 \label{eq:m''A<e}
 \mu  (A'')< \varepsilon
 \end{equation}
 and
 \begin{equation}
 \label{eq:bA''contClO}
 b^{\mu , h} (A'') = 
 \overline\Omega \cap \spt\mu.
 \end{equation}
 Regarding \eqref{eq:m''A<e}, we notice that it trivially follows from \eqref{eq:A'contQmenoO}, \eqref{eq:A'capOcontA'menoECC} and \eqref{eq:defA''}. As far as \eqref{eq:bA''contClO} is concerned, the inclusion $b^{\mu , h} (A'') \subset
 \overline\Omega \cap \spt\mu$  is immediately obtained as in \eqref{eq:bAsubsClOintS}.
 Moreover, since $b^{\mu , h} (A'')\supset b^{\mu , h} (A)\supset\Omega\cap\spt\mu$, we only need to show that
 \begin{equation} 
 \label{eq:bA''contdOcapS}
 b^{\mu , h} (A'')\supset
 \partial\Omega\cap\spt\mu
 \end{equation}
 to complete the proof of \eqref{eq:bA''contClO}. So let us consider $x\in  \partial\Omega\cap\spt\mu$ and observe that $\partial\Omega\subset A'$. Then $B_r(x)\subset A'$, provided $r$ is small enough, hence
 \begin{equation*}
 \label{eq:OB=A''B}
 \Omega\cap B_r(x)
 \supset
 A''\cap B_r(x)
 \supset
 A'\cap \Omega\cap B_r(x)
 =\Omega\cap B_r(x).
   \end{equation*}
   But we have also $x\in b^{\mu , h}(\Omega)$ 
   (by \eqref{eq:DOcapsptmsubb}) and thus   we obtain $x\in b^{\mu , h}(A'')$, which proves \eqref{eq:bA''contdOcapS}.
   
   \medskip
    \underline{Proof of (2).} 
    This statement is proved by the same argument used to prove (2) of \cite[Prop.3.2]{D32}, with some trivial adaptations.
 \end{proof}

 \begin{remark}
 Obviously, condition \eqref{eq:rpmuBrq} only makes sense if $q\leq p$. Moreover, if $q > n$ this condition implies that $\spt\mu$ is empty. In fact, if we assume $\spt\mu \not = \emptyset$ (and $q>n$), then we obtain the following contradiction:   
 \begin{itemize}
 \item On the one hand, as observed at the beginning of the proof of Theorem \ref{thm:GenProp3.2ofD32}, one would have $\mu (\spt\mu)>0$;
 \item On the other hand, by \cite[Th.6.9]{Mattila1995}, we have $\mu (\spt\mu)=0$.
 \end{itemize}
 These considerations make it clear why we assumed $q\leq \min\{n , p\}$ in Theorem \ref{thm:GenProp3.2ofD32}.  
 \end{remark}

 \begin{remark}
 Let $p,q$ be as in Theorem \ref{thm:GenProp3.2ofD32}. Then it is easy to verify that $\frac{np}{q}-q=0$ if and only if $p=q=n$. 
 \end{remark}

 \begin{remark}
 \label{rem:Cond4.3inClassCase}
 We observe that:
 \begin{itemize}
 \item[(1)] If $\mu=\leb^n$ then condition \eqref{eq:DOcapsptmsubb} is verified whenever $\partial \Omega$ is Lipschitz (for all $h\in [0,+\infty)$). Hence Theorem \ref{thm:GenProp3.2ofD32} yields immediately \cite[Prop.3.2]{D32}.
 \item[(2)] No regularity assumption on $\partial\Omega$ will suffice to ensure that condition \eqref{eq:DOcapsptmsubb} is verified for all $\mu\in\R$. For example, if $\Omega$ is a ball and $\mu:=\hf^{n-1}\ree\partial\Omega$ then $\partial\Omega\cap\spt\mu=\partial\Omega$ and $b^{\mu , h}(\Omega)=\emptyset$ (for all $h\in [0,+\infty)$).
 \end{itemize}
 \end{remark}

 \begin{remark}
 Let $\mu:=\hf^k\ree S$, where $S$ is an open imbedded $k$-submanifold of $\rn^n$ of class $C^1$ with $k\leq n-1$. Moreover let $\partial \Omega$ be of class $C^1$ and assume that $S$ and $\partial\Omega$ meet transversely at $x$, namely 
 \begin{equation*}
 x\in \partial\Omega\cap S,\quad
 \dim (T_xS + T_x(\partial\Omega))=n,
 \end{equation*}
 where $T_xS$ and $T_x(\partial\Omega))$ are the tangent space of $S$ at $x$ and the tangent space of $\partial\Omega$ at $x$, respectively. We observe that then we also have $\dim (T_xS\cap T_x(\partial\Omega))=k-1$ and this fact implies that near $x$ the set $\partial\Omega\cap S$ is an imbedded $(k-1)$-submanifold of $\rn^n$ of class $C^1$. Then, with a standard argument based on the area formula, we can prove that $x\in b^{\mu , 0}(\Omega)$ (hence $x\in b^{\mu , h}(\Omega)$ for all $h\in [0,+\infty)$). Therefore, if we now assume that $S$ and $\partial\Omega$ meet transversely everywhere (i.e., at every point in $\partial\Omega\cap S$), then we find $\partial\Omega\cap S\subset b^{\mu , 0}(\Omega)$. This does not imply that condition \eqref{eq:DOcapsptmsubb} is verified. For example, consider the case $n:=3$, $k:=2$ and
 \begin{equation*}
 \Omega := B_1(0),\quad
 S:=\{
 (x_1,x_2,x_3)\in\rn^3\,\vert\,  x_3(x_3-1)=0
 \}
 \setminus \{(0,0,1)\}.
 \end{equation*}
 In this case $S$ and $\partial\Omega$ meet transversely everywhere and $\partial\Omega\cap S = b^{\mu , h}(\Omega)$ for all $h\in [0,+\infty)$. Hence we have also
 $$
 (0,0,1)\not\in b^{\mu , h}(\Omega)
 $$
 and
 \begin{equation*}
 \partial\Omega \cap\spt\mu
 = \partial\Omega \cap\overline S
 = ( \partial\Omega \cap S) \cup \{ (0,0,1)\}
 = b^{\mu , h}(\Omega) \cup \{ (0,0,1)\}
 \end{equation*}
 for all $h\in [0,+\infty)$.
  \end{remark}

 \begin{remark}
 It is natural to ask whether Theorem \ref{thm:GenProp3.2ofD32} can be extended to the case that $\bar r$ depends on $x\in\spt\mu$. After trying to prove such a generalisation, we are inclined to believe that the answer is negative, but we have no counterexamples. 
 \end{remark}
 

 \section{Applications of Theorem \ref{thm:GenProp3.2ofD32}, two remarkable examples}

  \subsection{First example}
 Let $\mu=\leb^n$ and $p=q=n$. Then, applying Theorem \ref{thm:GenProp3.2ofD32} and recalling (1) of Remark \ref{rem:Cond4.3inClassCase}, we obtain \cite[Prop.3.2]{D32}.

 \subsection{Second example}
 \label{subsect:2ndApplsOfFirstMainThm}
 Let $k\leq n$ and consider a bounded open subset $G$ of $\rn^k$, with boundary of class $C^1$. Let $\varphi\in C^1 (\rn^k , \rn^n)$ be such that $\varphi\vert_G$ is injective and
 \begin{equation*}
 J\varphi (y)=\left[
 \det\left[
 (D\varphi)^t\times (D\varphi)
 \right] (y)
 \right]^{1/2}
 >0
 \end{equation*}
 for all $y\in \overline G$. We will apply Theorem \ref{thm:GenProp3.2ofD32} to the measure $\mu:=\hf^k\ree\varphi (G) = \hf^k\ree\varphi (\overline G)$, but in order to do so, we must first prove the following result.
 \begin{proposition}
 \label{prop:rk/CleqmBleqCrk}
 There exist $C, \bar r\in (0, +\infty)$ such that
 \begin{equation}
 \label{eq:rk/CleqmBleqCrk}
 \frac{r^k}{C}
 \leq
 \mu (B_r(x))
 \leq
 C r^k
 \end{equation}
 for all $x\in\spt\mu$ and $r\in (0,\bar r]$.
 \end{proposition}
 \begin{proof}
 Let us first consider $y\in \overline G$ and observe that the number
 \begin{equation*} 
 \lambda (y):= \min \left\{\left([(D\varphi)^t\times (D\varphi)] (y)u\right)\cdot u\,\vert  \, u\in\rn^k,\, \vert u\vert=1\right\}
 \end{equation*}
 is the smallest eigenvalue of the matrix
 \begin{equation*}
 [(D\varphi)^t\times (D\varphi)] (y)
 \end{equation*}
 cf.  \cite[Ch.VIII, Th.3.3]{Lang1987}. Hence and recalling that the zeros of a monic polynomial depend continuously on its coefficients (cf. \cite[Sect.1.3]{RahmanShmeisser2002}) we obtain that the function $\lambda:\overline G \to \rn$ is continuous. Then also the function mapping $y\in\overline G$ to
 \begin{equation*}
 m_0 (y):= \frac{1}{2}\left[\min \left\{\left([(D\varphi)^t\times (D\varphi)] (y)u\right)\cdot u\,\vert  \, u\in\rn^k,\, \vert u\vert=1\right\}\right]^{1/2}
 \end{equation*}
 has to be continuous. Since $\overline G$ is compact, there exists $y_0\in\overline G$ such that
 \begin{equation*}
 m_{00}:= m_0(y_0)= \min_{y\in\overline G} m_0(y).
 \end{equation*}
 Observe that $m_{00}>0$ by Lemma \ref{lem:minLuu>0}. Furthermore, since $D\varphi$ is continuous, we easily see that there must exist $r_0\in (0,1]$ such that
 \begin{equation*}
\sigma_r (y):= \max_{z\in\overline{B_{r}(y)}}\Vert (D\varphi) (z) - (D\varphi) (y)\Vert\leq m_{00}, 
 \end{equation*}
for all $y\in\overline G$ and $r\in (0,r_0]$, where $\Vert \cdot \Vert$ denotes the Hilbert-Schmidt norm of matrices. Now using inequality \eqref{eq:deltaphigeqmr-sr} we obtain
 \begin{equation}
 \label{eq:CorrTo(3.8)}
 \vert\varphi (z) - \varphi (y)\vert
 \geq 2m_0(y) r - \sigma_r (y) r
 \geq m_{00}r,
 \end{equation}
 for all $y\in\overline G$, $z\in\partial B_r(y)$ and $r\in (0,r_0]$. On the other hand, recalling \eqref{eq:deltaphileqMr}, we also have
 \begin{equation}
 \label{eq:CorrTo(3.4)}
 \vert\varphi (z) - \varphi (y)\vert
 \leq m_1 r
 \end{equation}
 for all $y\in\overline G$, $z\in\partial B_r(y)$ and $r\in (0,1]$, where $m_1$ is defined as in \eqref{def:Kandm1}.  From \eqref{eq:CorrTo(3.8)} and \eqref{eq:CorrTo(3.4)} it follows that
 \begin{equation}
 \label{eq:FundamEstim}
 \varphi (\overline G)\cap B_{m_{_{00}}r}(\varphi (y))
 \subset
 \varphi (\overline G)\cap\varphi (B_r(y))
 \subset
  \varphi (\overline G)\cap B_{m_{_{1}}r}(\varphi (y)),
 \end{equation}
 for all $y\in\overline G$ and $r\in (0,r_0]$. Now, using \eqref{eq:FundamEstim} , we can proceed to the proof of \eqref{eq:rk/CleqmBleqCrk}:
 \vskip2mm
 \begin{itemize}
 \item We first prove by contradiction the following claim: There exist  $C_1, r_1\in (0,+\infty)$ such that
\begin{equation}
 \label{eq:rk/CleqmBleqCrk_1}
  \mu (B_r(x))
  \geq
 \frac{r^k}{C_1}
 \end{equation}
 for all $x\in\spt\mu =\varphi (\overline G)$ and $r\in (0, r_1]$. If this were not true, for each positive integer $j$ there would exist $y_j\in \overline G$ and $\rho_j\in (0, 1/j]$ such that
 \begin{equation}
 \label{eq:HkfGB<rk/j}
 \hf^k\left(
 \varphi (G)\cap B_{\rho_j}(\varphi (y_j))
 \right)
 <
  \frac{\rho_j^k}{j}.
 \end{equation}
 Since $\overline G$ is compact we can assume that $y_j\to \bar y\in\overline G$, as $j\to +\infty$. On the other hand, by the second inclusion in \eqref{eq:FundamEstim} and the area formula, we have
 \begin{equation*}
  \hf^k\left(
 \varphi (G)\cap B_{\rho_j}(\varphi (y_j))
 \right)
 \geq
 \hf^k\left(
\varphi (G)\cap\varphi (B_{\rho_j/m_1}(y_j)
 \right)
 =
 \int_{G \cap B_{\rho_j/m_{_{00}}}(y_j)}J\varphi\, d\leb^k,
 \end{equation*}
 provided $j$ is large enough. Hence, recalling that $\partial G$ is of class $C^1$, we find
 \begin{equation*}
 \liminf_{j\to +\infty}\frac{ \hf^k\left(
 \varphi (G)\cap B_{\rho_j}(\varphi (y_j))
 \right)}{\leb^k (B_{\rho_j/m_{_{00}}}(y_j))}
 \geq
 \frac{J\varphi (\bar y)}{2}>0
 \end{equation*}
 which contradicts \eqref{eq:HkfGB<rk/j}. Thus the claim above has to be true.
 
 \vskip2mm
 \item From the first inclusion in \eqref{eq:FundamEstim} and the area formula it follows that
 \begin{equation*}
 \hf^k\left(
 \varphi (G)\cap B_r(\varphi (y))
  \right)
 \leq
 \hf^k \left(
 \varphi (G)\cap \varphi (B_{r/m_{_{00}}}(y))
  \right)
  =\int_{G \cap B_{r/m_{_{00}}}(y)}J\varphi\, d\leb^k
 \end{equation*}
 for all $y\in\overline G$ and $r\in (0, m_{_{00}}r_0]$.  Thus, since $J\varphi$ is bounded in $\overline G$, there must exist a positive constant $C_2$ (which does not depend on $x$ and $r$) such that
 \begin{equation}
  \label{eq:rk/CleqmBleqCrk_2}
 \mu (B_r(x))
 =
 \hf^k\left(
 \varphi (G)\cap B_r(x)
  \right)
  \leq
  C_2 r^k
 \end{equation}
  for all $x\in \spt\mu=\varphi(\overline G)$ and $r\in (0, m_{_{00}}r_0]$.
  
  \vskip2mm
 \item Finally, the inequalities \eqref{eq:rk/CleqmBleqCrk_1} and \eqref{eq:rk/CleqmBleqCrk_2} yield \eqref{eq:rk/CleqmBleqCrk} with $C:=\max\{C_1,C_2\}$ and $\bar r :=\min\{r_1, m_{_{00}}r_0\}$.
 \end{itemize}
 \end{proof}
 
  Now, by applying Theorem \ref{thm:GenProp3.2ofD32} with $\mu=\hf^k\ree\varphi (G)$ and $p=q=k$ (taking Proposition \ref{prop:rk/CleqmBleqCrk} into account), we obtain:
  
  \begin{corollary}
  The following properties hold  for all $\varepsilon >0$ and $h> n-k$: 
 \begin{itemize}
 \item[(1)] If $\Omega\subset\rn^n$ is a bounded open set, then there exists an open set $A\subset\Omega$ such that
 \begin{equation*}
 \label{eq:(3.2)inD32_cor}
 \hf^k(\varphi (G) \cap A)< \varepsilon,
 \qquad
 \Omega\cap \varphi(\overline G) \subset 
 b^{\hf^k\ree\varphi (G), h} (A)
 \subset
\overline\Omega\cap \varphi(\overline G).
 \end{equation*}
 In the special case when
 \begin{equation*}
 \label{eq:DOcapsptmsubb_cor}
 \partial\Omega\cap\varphi(\overline G)\subset b^{\hf^k\ree\varphi (G),h}(\Omega),
 \end{equation*}
 the set $A$ can be chosen so that we have
 \begin{equation*}
 \label{eq:(3.2)inD32bis_cor}
  b^{\hf^k\ree\varphi (G), h} (A)= \overline\Omega\cap \varphi(\overline G).
 \end{equation*}
 \item[(2)] There is an open set $U\subset\rn^n$ satisfying
 $$
 \hf^k (\varphi (G)\cap U)<\varepsilon,\qquad
 b^{\hf^k\ree\varphi (G),h}(U) 
 =\varphi(\overline G).
 $$
 \end{itemize}
  \end{corollary}

 \section{Density degree functions}
 
 Let $\mu\in\R$ be non-trivial, i.e., $\spt\mu\not = \emptyset$. We will follow the path traced in \cite{D34}.
 
 \medskip
 First of all observe that if 
 $E\subset \rn^n$ and $x\in\rn^n$,
 then the set $\{h\in [0,+\infty)\,\vert\, x\in E^{\mu , h}\}$ is a (possibly empty) interval.

\begin{definition}
\label{def:IntClBound}
Let $E$ be a subset of $\rn^n$. Then the {\it density degree of $E$ (w.r.t. $\mu$)} is the function $d_E^\mu:\rn^n\to \{-n\}\cup [0,+\infty]$  defined as follows:
$$
d_E^\mu (x):=
\begin{cases}
        \sup\{h\in [0,+\infty)\,\vert\, x\in E^{\mu , h}\}&  \text{if  $x\in E^{\mu , 0}$}\\
        -n& \text{if  $x\not\in E^{\mu , 0}$}.
\end{cases}
$$
For $m\in [0,+\infty]$ we also define
\begin{equation*}
\essint^{\mu , m}E := \left\{
x\in\rn^n\,\vert\, d_E^\mu(x) >m
\right\},\qquad
\esscl^{\mu , m}E := \left\{
x\in\rn^n\,\vert\, d_E^\mu(x) \geq m
\right\}
\end{equation*}
and
$$
\partial^{\mu , m}E:= \esscl^{\mu , m}E\setminus \essint^{\mu , m}E
= \left\{
x\in\rn^n\,\vert\, d_E^\mu (x) = m
\right\}.
$$
 When the following identity holds
$$
E^{\mu , 0}\equalmm \partial^{\mu , m}E= \left\{
x\in\rn^n\,\vert\, d_E^\mu (x) = m
\right\}
$$
we say that $E$ is a {\it uniformly $(\mu, m)$-dense set}.
\end{definition}

\begin{remark}
\label{rem:TFaboutdE}
The following trivial facts occur:
\begin{itemize}
\item[(1)] If $E\equalmm \emptyset$ then $E^{\mu , 0}=\emptyset$ and hence $d_E^\mu\equiv -n$;
\item[(2)] $\esscl^{\mu , 0}E = E^{\mu , 0}$;
\item[(3)] $\essint^{\mu , +\infty}E=\emptyset$, hence $\partial^{\mu , +\infty}E:= \esscl^{\mu , +\infty}E$.
\end{itemize}
\end{remark}

\begin{example}
If $E$ is open, then $d_E^\mu (x)=+\infty$ for all $x\in E$. Hence
$$
E\subset \essint^{\mu , m} E
$$
for all $m\in [0,+\infty)$. Observe that the strict inclusion can occur, e.g., for $\mu:=\leb^n$ and $E:= B_r\setminus \{0\}$ (in such a case one has $ \essint^{\mu , m} E=B_r$).
\end{example}


This proposition collects some very simple (nevertheless interesting) facts.

\begin{proposition}
\label{prop:FirstSimpleProp}
Let $E$ be a subset of $\rn^n$ and $m\in [0,+\infty]$. The following properties hold:
\begin{itemize}
\item[(1)] $\partial^{\mu , k}E\cap \partial^{\mu , m}E=\emptyset$, if $k\in [0,+\infty]$ and $k\not = m$.
\vskip2mm
\item[(2)] $\essint^{\mu , m} E =\bigcup_{k>m} E^{\mu , k}$.
\vskip2mm
\item[(3)] If $m>0$ then $\esscl^{\mu , m} E = \bigcap_{l\in [0,m)}E^{\mu , l}$.
\vskip2mm
\item[(4)] $\essint^{\mu , m} E$, $\esscl^{\mu , m} E$ and $\partial^{\mu , m}E$ are $\mu$-measurable sets.
\vskip2mm
\item[(5)] $\essint^{\mu , m}E\subset E^{\mu , m}\subset \esscl^{\mu , m} E$.
\vskip2mm
\item[(6)] The following two claims are equivalent:
\begin{itemize}
\item[$\bullet$] $E$ is a uniformly $(\mu ,m)$-dense set;
\item[$\bullet$] $\esscl^{\mu , m}E\equalmm E^{\mu , 0}$ and $\essint^{\mu , m}E\equalmm \emptyset$.
\end{itemize}
\item[(7)]   $E$ is a uniformly $(\mu , 0)$-dense set if and only if $\essint^{\mu , 0}E\equalmm\emptyset$.   
\vskip2mm
\item[(8)] The function $d_E^\mu$ is measurable.
\end{itemize}
\end{proposition}

\begin{proof}
Definition \ref{def:IntClBound} yields at once (1), (2) and (3). Statement (4) follows trivially from (2) and (3), by recalling (2) in Remark \ref{rem:TFaboutdE}, (2) in Remark \ref{rem:PropSuperd} and  (1) in Proposition \ref{prop:ProprBaseTop}. Also (5) follows trivially from (2) and (3), by recalling (2) in Remark \ref{rem:PropSuperd}.

\medskip
 Let us prove (6). 
 \begin{itemize}
 \item If we assume that the first claim is true, then, by recalling also (3), we obtain
 $$
 \esscl^{\mu , m}E\subset E^{\mu , 0} \equalmm    \esscl^{\mu , m}E\setminus \essint^{\mu , m} E\subset \esscl^{\mu , m}E.
 $$
 This proves the first formula in the second claim. It also proves that $\esscl^{\mu , m}E\equalmm \esscl^{\mu , m}E\setminus \essint^{\mu , m}$, hence the last formula in the second claim follows by recalling (5).
 \item Conversely, if we assume that the second claim is true, then
 $$
 \partial^{\mu , m}E =  \esscl^{\mu , m}E\setminus \essint^{\mu , m}E\equalmm E^{\mu , 0}
 $$
 i.e. $E$ is a uniformly $(\mu ,m)$-dense set.
\end{itemize}
 
 \medskip
Now, the statement (7) follows at once from (2) in Remark \ref{rem:TFaboutdE} and (6).
Finally, observe that for $a\in\rn$ one has
$$
\left\{
x\in\rn^n\,\vert\, d_E^\mu(x)\geq a
\right\}
=
\begin{cases}
       \rn^n& \text{if  $a\leq -n$}\\
        E^{\mu , 0}&  \text{if $a\in (-n,0)$}\\
        \esscl^{\mu , a} E &  \text{if $a \geq 0$}
\end{cases}
$$
by Definition \ref{def:IntClBound}.  Hence (8) follows from (1) in Proposition \ref{prop:ProprBaseTop} and (4).
\end{proof}

\begin{remark}
Let $E\subset\rn^n$ and $m\in [0,+\infty ]$. Then, from (4) and (6) of Proposition \ref{prop:FirstSimpleProp} and (5) in Remark \ref{rem:PropSuperd} it follows that the following statements are equivalent:
\begin{itemize}
\item $E\equalmm \partial^{\mu , m}E$;
\item $E\in\M_\mu$ and $E$ is a uniformly $(\mu ,m)$-dense set,
\item $ \esscl^{\mu , m} E\equalmm E$ and $\essint^{\mu , m}E\equalmm \emptyset$.
\end{itemize}
\end{remark}

\begin{remark}
Proposition \ref{prop:FirstSimpleProp} holds whatever negative value is assigned, in Definition \ref{def:IntClBound}, to the restriction of $d_E^\mu$ to $E^{\mu , 0}$. We chose $-n$ only because this way the function $n+d_E^{\leb^n}$ coincides with the density degree function $d_E$ defined in \cite[Def.5.1]{D34}.
\end{remark}

The following proposition is an easy consequence of (1) and (4) of Proposition \ref{prop:FirstSimpleProp} (cf. \cite[Prop.5.2]{D34}).

\begin{proposition}
Let $E$ be a measurable subset of $\rn^n$. Then the set
$$
\left\{
m\in [0,+\infty ]\,\vert\, \mu(\partial^{\mu , m}E)>0
\right\}
$$
is at most countable.
\end{proposition}

Now we prove a result about approximation of a set, given as the closure of an open set, by closed subsets having small density degree (w.r.t. $\mu$). The proof is obtained by adapting the argument used in \cite[Prop.5.4]{D34}.

\begin{proposition}
\label{cor:ThinSetF}
Assume that:
\begin{itemize}
\item[(i)] There exist $C, p, q, \bar r \in (0, +\infty)$ such that $q\leq \min\{n , p\}$ and
\begin{equation*}
\frac{r^p}{C}
\leq  \mu (B_r(x))
\leq Cr^q
\end{equation*}
for all $x\in\spt \mu$ and $r\in (0,\bar r)$;
\item[(ii)] It is given a non-empty bounded open set $\Omega\subset \rn^n$ with the following property: There exists an open bounded set $\Omega'\subset\rn^n$ such that $\Omega\subset\Omega'$ and $\partial\Omega' \cap \spt\mu\subset b^{\mu , h} (\Omega')$ for all $h>\overline m:= \frac{np}{q} - q$.
\end{itemize}
Then for all $H\in \left(0, \mu \left(\overline\Omega\right)\right)$ there exists a closed   subset $F$ of $\overline\Omega$ such that 
$$
\mu (F) > H, \quad
\essint^{\mu , \overline m}F\equalmm \emptyset.
$$
\end{proposition}

\begin{proof}
Let $j$ be an arbitrary positive integer. Then, by Theorem \ref{thm:GenProp3.2ofD32}, there exists an open set $A_j\subset \Omega'$ such that
 \begin{equation}
 \label{eq:estmAjebmhj}
 \mu(A_j)
 <
 {\mu \left(\overline\Omega\right) - H \over 2^j } , \qquad  b^{\mu , h_j}(A_j)= \overline{\Omega'}\cap\spt\mu
 \end{equation}
 with
 $$
 h_j:=\overline m +\frac{1}{j}=\frac{np}{q} - q+\frac{1}{j}.
 $$
 Define
 $$
 K_j:=\overline{\Omega'}\cap A_j^c,\qquad K:=\bigcap_{j=1}^{+\infty}K_j=\overline{\Omega'}\cap\left(\bigcup_{j=1}^{+\infty} A_j\right)^c.
 $$
 Then $K$ is closed and
 \begin{equation}
 \label{eq:mK>Dm+H}
 \begin{split}
 \mu(K)
 &=  \mu \left(\overline{\Omega'}\right) -
  \mu\left(\cup_j A_j\right)
 \geq  
 \mu \left(\overline{\Omega'}\right) - \sum_j  \mu(A_j)
 >   
 \mu \left(\overline{\Omega'}\right) -  \mu\left(\overline\Omega\right) + H,
 \end{split}
 \end{equation}
 by \eqref{eq:estmAjebmhj}.
 Moreover, by (2), (3), (5) of Proposition \ref{prop:mainPropOfbmuh} and \eqref{eq:estmAjebmhj}, we have
 \begin{equation*}
 \begin{split}
 K_j^{\mu , h_j} 
 &\equalmm 
 \left[ b^{\mu , h_j}(K_j^c)\right]^c
 =
 \left[ b^{\mu , h_j}\left( \left(\overline{\Omega'}\right)^c \cup A_j \right)\right]^c\\
  &=
  \left[ 
  b^{\mu , h_j}\left(\left(\overline{\Omega'}\right)^c\right)\cup b^{\mu , h_j}( A_j) \right]^c
  \\
  &\subset
  \left[\left( \left(\overline{\Omega'}\right)^c\cap\spt\mu\right)
  \cup
  \left(\overline{\Omega'}\cap\spt\mu\right)\right]^c\\
  &=
  (\spt \mu)^c
 \end{split}
 \end{equation*}
 that is
 \begin{equation*}
 \label{eq:Fjmhj=empty}
 K_j^{\mu , h_j} 
 \equalmm 
 \emptyset
 \end{equation*}
 for all $j$. Moreover, for each $k\in (\overline m, +\infty)$ we can find $j$ such that $k> h_j$, hence
 \begin{equation*}
 \label{eq:FmksubFjmhj}
 K^{\mu , k}\subset  K^{\mu , h_j}
 \subset
  K_j^{\mu , h_j}\equalmm\emptyset,
 \end{equation*}
by (2) of Remark \ref{rem:PropSuperd}.
Recalling (2) of Proposition \ref{prop:FirstSimpleProp}, we obtain
 $$
 \essint^{\mu , \overline m}K=\bigcup_{k>\overline m}K^{\mu , k}\equalmm \emptyset.
 $$
 Now define
 $$
 F:=\overline\Omega\cap K.
 $$
 Then $F$ is a closed subset of $\overline\Omega$ and (again by (2) of Proposition \ref{prop:FirstSimpleProp})
 \begin{equation*}
 \essint^{\mu , \overline m} F
 \subset
 \essint^{\mu , \overline m} K
 \equalmm 
 \emptyset, \text{ i.e. } \essint^{\mu , \overline m} F \equalmm \emptyset.
 \end{equation*}
 Moreover
 \begin{equation*}
 \mu (F) 
 = 
 \mu(K)-\mu\left( K\setminus\overline\Omega\right)
 >
  \mu \left(\overline{\Omega'}\right) -  \mu\left(\overline\Omega\right) + H -\mu\left( K\setminus\overline\Omega\right)
 \end{equation*}
 by \eqref{eq:mK>Dm+H}, where
 \begin{equation*}
 \mu \left(\overline{\Omega'}\right) -  \mu\left(\overline\Omega\right)
 =
 \mu\left(\overline{\Omega'}\setminus\overline\Omega\right)
 \geq
 \mu \left(K\setminus\overline\Omega\right).
 \end{equation*}
 Hence $\mu (F)>H$.
  \end{proof}


\begin{remark}
\label{rem:Prop7.3ImplProp5.4old}
From Proposition \ref{cor:ThinSetF} and (1) of Remark \ref{rem:Cond4.3inClassCase}, we easily obtain  \cite[Prop.5.4]{D34}.
\end{remark}

\begin{corollary}
\label{cor:intl0E=empty}
Let $G$ and $\varphi$ be as in Section \ref{sect:ExampleC1Map}. Moreover, let $A\subset \varphi (G)$ be open with respect to the topology induced in $\varphi (G)$ by $\tau (\rn^n)$ and assume that 
\begin{equation}
\label{eq:clAsubphiG}
\overline A\subset\varphi (G).
\end{equation}
Then for all $H\in \left(0,\hf^k \left(\overline A\right)\right)$ there exists a closed set $E\subset \overline A$ such that
\begin{equation}
\label{eq:lF>HandintF=empty}
\lambda (E) = \hf^k (E) > H,
\quad
\essint^{\lambda , 0}E
= \hskip-3.4mm^{^\lambda}
\; \emptyset
\end{equation}
where $\lambda:=\hf^k\ree\varphi (G)$. In particular, $E$ is a uniformly $(\lambda , 0)$-dense set.
\end{corollary}

\begin{proof}
Let us consider the bounded open set
\begin{equation*}
D :=(\varphi\vert_G)^{-1}(A)
\end{equation*}
and observe that, by \eqref{eq:clAsubphiG}, we have also
\begin{equation}
\label{eq:clPA=PclA}
\overline D
=
\overline{(\varphi\vert_G)^{-1}(A)}
=
 (\varphi\vert_G)^{-1}(\overline A)\subset G.
\end{equation}
Now let $H\in \left(0,\hf^k \left(\overline A\right)\right) =  \left(0,\lambda \left(\overline A\right)\right)$ and consider $H'\in \left(0, \leb^k \left(\overline D\right)\right)$ satisfying
\begin{equation}
\label{eq:H'geqLD-quot}
H' \geq \leb^k\left(
\overline D\right)
-
\frac{\lambda \left(\overline A\right) - H}{M},
\end{equation}
where 
$$
M:=\max_{\overline D} J\varphi.
$$
From Proposition \ref{cor:ThinSetF} (with $n=p=q=k$ and $\mu=\leb^k$) and recalling (1) of Remark \ref{rem:Cond4.3inClassCase}, it follows that a closed set $K\subset\overline D$ has to exist such that
\begin{equation}
\label{eq:LkK>H'andint=empty}
\leb^k (K) 
> H', 
\quad
\essint^{\leb^k , 0}K
= \hskip-3.9mm^{^{\leb^k}}
\; \emptyset.
\end{equation}
Then consider $h\in [0,+\infty)$ and the closed set
\begin{equation*}
E:=\varphi (K).
\end{equation*}
Observe that
\begin{equation}
\label{eq:FllsubFsubPG}
E^{\lambda , h}
\subset
E
\subset
\varphi \left( \overline D\right)
=
\overline A
\subset \varphi (G),
\end{equation}
by \eqref{eq:clPA=PclA} and (7) in Remark \ref{rem:PropSuperd}. Hence and by the area formula (cf. \cite[Cor. 5.1.13]{KrantzParks2008}) we obtain
\begin{equation}
\label{eq:lclA-lF}
\begin{split}
\lambda\left(\overline A\right)
-
\lambda (E)  
&=
\hf^k\left(\varphi\left(\overline D\right)\right)
-
\hf^k (\varphi (K))=\int_{\overline D\setminus K} J\varphi\, d\leb^k\\
&\leq  
M
\left(
\leb^k\left(\overline D\right)
-
\leb^k (K)
\right).
\end{split}
\end{equation}
The inequality in
\eqref{eq:lF>HandintF=empty}
now follows easily from \eqref{eq:H'geqLD-quot}, \eqref{eq:LkK>H'andint=empty} and \eqref{eq:lclA-lF}. 

\medskip
From Proposition \ref{prop:phi(E(k+h))=E(Hkphih)}, Remark \ref{rem:P-1Ek+h=PG-1Ek+h}, \eqref{eq:LkK>H'andint=empty}, (2) in Proposition \ref{prop:FirstSimpleProp} and \eqref{eq:FllsubFsubPG},  also taking into account (1) and (2) in Remark \ref{rem:PropSuperd}, it follows that  
\begin{equation*}
E^{\lambda , h}
=
E^{\lambda , h}\cap \varphi (G)
=\varphi\left(
[(\varphi\vert_G)^{-1}(E)]^{(k+h)}\cap G
\right)
= \varphi\left(
K^{(k+h)}\cap G
\right)
{= \hskip-3.4mm^{^\lambda}}\;
\emptyset
\end{equation*}
for all $h\in (0, +\infty)$. Hence, recalling again (2) in Proposition \ref{prop:FirstSimpleProp} and (2) in Remark \ref{rem:PropSuperd}, we obtain 
$$
\essint^{\lambda , 0}E
= \hskip-3.4mm^{^\lambda}
\; \emptyset.
$$ 
Finally, $E$ is a uniformly $(\lambda , 0)$-dense set, by (7) of Proposition \ref{prop:FirstSimpleProp}.
\end{proof}

  \begin{remark}
  In general, Proposition \ref{cor:ThinSetF} does not provide the optimal result. For example, if we apply Proposition \ref{cor:ThinSetF} directly to the measure $\lambda$ carried by a $k$-dimensional imbedded $C^1$ submanifold of $\rn^n$ with $C^1$ boundary we get a worse result than that obtained in Corollary \ref{cor:intl0E=empty}. To verify this fact, let us consider $G$ and $\varphi$ as in Section \ref{sect:ExampleC1Map} and further assume that $\partial G$ is of class $C^1$.We observe that hypothesis (i) of Proposition \ref{cor:ThinSetF} is verified, with 
  $$
  \mu =\lambda =\hf^k\ree\varphi (G), \quad
  p=q=k,
  $$
  by Proposition \ref{prop:rk/CleqmBleqCrk}. 
  Now let $A\subset \varphi (G)$ be open with respect to the topology induced in $\varphi (G)$ by $\tau (\rn^n)$ and assume that \eqref{eq:clAsubphiG} holds. By a standard argument, it follows that an open set $\Omega\subset\rn^n$ exists such that 
  \begin{equation*}
  A=\Omega \cap \varphi (G),
  \quad 
  \overline A=\overline\Omega \cap \varphi (G).
  \end{equation*}
  Since $\spt\mu$ is bounded, there is an open ball $B\subset\rn^n$ such that $\Omega\subset B$ and $\partial B\cap\spt\mu=\emptyset$. Hence (ii) of Proposition \ref{cor:ThinSetF} is trivially verified, with $\Omega' = B$. Now consider any $H\in \left( 0 , \hf^k\left(\overline A\right)\right)$ and observe that $\hf^k\left(\overline A\right)=\lambda \left(\overline \Omega\right)$. Then, by Proposition \ref{cor:ThinSetF}, there exists a closed   subset $F$ of $\overline\Omega$ such that 
$\lambda (F) > H$ and $\essint^{\lambda , n-k}F = \hskip-3.4mm^{^\lambda} \; \emptyset$, i.e.,
$$
\lambda (E) > H, \quad
\essint^{\lambda , n-k}E = \hskip-3.4mm^{^\lambda}
\; \emptyset
$$
where $E:=F\cap \varphi (G)$, which is closed with respect to the topology induced in $\varphi (G)$ by $\tau (\rn^n)$. Therefore, this argument does not prove the result obtained in Corollary \ref{cor:intl0E=empty}, namely, that there are closed subsets of $\overline A$ of arbitrarily close measure to $\hf^k\left(\overline A\right)$ that are also uniformly $(\lambda , 0)$-dense.
   \end{remark}

    \section{A Schwarz-type result}
  
  We will prove the following result that generalizes the classical Schwarz theorem on cross derivatives (cf. Remark \ref{rem:ThmReducesToSchwarz} below).

  \begin{theorem}
  \label{thm:gengenschwarz}
  Let us consider $\mu\in\R$, an open set $\Omega\subset\rn^n$, $f, G, H\in C^1(\Omega)$, a couple of integers $p,q$ such that $1\leq p< q\leq n$ and $x\in\rn^n$. Assume that:
  \begin{itemize}
  \item[(i)] For $i=p, q$, the $i$-th distributional derivative of $\mu$ is a Borel real measure on $\rn^n$ also denoted $D_i\mu$ (with no risk of misinterpretation), so that we have $D_i\mu (\varphi)=-\int D_i\varphi\, d\mu = \int\varphi\, d(D_i\mu)$, for all $\varphi\in C_c^1(\rn^n)$;
  \vskip2mm
  \item[(ii)] $x\in\Omega\cap A^{\mu , 1}$, where $A:=\{ y\in\Omega\,\vert\, (D_pf (y), D_qf(y)) =(G(y),H(y))\}$ (in particular $x\in\spt\mu$);
   \vskip2mm
  \item[(iii)] $\lim_{\rho\to 1-}\sigma (\rho)=1$, where $\sigma (\rho):= \liminf_{r\to 0+}\frac{\mu (B_r(x))}{\mu(B_{\rho r}(x))}$ (note that $\sigma$ is decreasing);
   \vskip2mm
  \item[(iv)] For $i=p,q$, one has $\lim_{r\to 0+}\frac{\vert D_i\mu\vert (B_r(x))}{r\mu (B_r(x))}=0$. 
  \end{itemize}
  Then $D_pH (x) = D_qG(x)$.
  \end{theorem}
  
  \begin{proof}
Let $\rho\in (0,1)$ and consider $g\in C_c^2(B_1(0))$ such that $0\leq g\leq 1$, $g\vert_{B_\rho(0)}\equiv 1$
and
$$
\vert D_ig\vert
\leq \frac{2}{1-\rho}\qquad (i=1,\ldots , n).
$$
For every real number $r$ such that $0<r<{\rm dist}(x,\rn^n\setminus \Omega)$, we define $g_r\in C_c^2(B_r(x))$ as
$$
g_r(y):=g\left(
\frac{y-x}{r}
\right),\quad y\in \rn^n
$$
and observe that (for all $y\in B_r(x)$ and $i=1,\ldots,n$)
\begin{equation}
\label{eq:EstGradDgr}
\vert D_i g_r(y)\vert
=
\frac{1}{r}\left\vert
 D_ig\left(
\frac{y-x}{r}
\right)
\right\vert
\leq \frac{2}{r(1-\rho)}.
\end{equation}
Moreover define 
$$
\Gamma:= D_pH - D_qG.
$$
Then, after a simple computation in which we use only (i), the definition of $A$ in (ii) and the identity $D_pD_q g_r = D_qD_p g_r$, we arrive at the following equality (where $B_r$ and $B_{\rho r}$ stand for $B_r(x)$ and $B_{\rho r}(x)$, respectively):
\begin{equation*}
\begin{split}
\int_{B_r}\Gamma g_r\, d\mu
&=
\int_{B_r} (g_r G + f D_p g_r)\, d (D_q\mu)
- \int_{B_r} (g_r H + f D_q g_r)\, d (D_p\mu)\\
&-\int_{B_r\setminus A}  (H - D_qf) D_p g_r \, d\mu
+ \int_{B_r\setminus A}  (G - D_pf) D_q g_r \, d\mu.
\end{split}
 \end{equation*}
 Hence, by also recalling the polar decomposition theorem (cf. \cite[Cor.1.29]{AFP2006}) and \eqref{eq:EstGradDgr}, we obtain
 \begin{equation*}
\begin{split}
\left\vert\int_{B_r}\Gamma g_r\, d\mu\right\vert
&\leq
\int_{B_r} \left(g_r \vert G\vert + \vert f\vert\, \vert D_p g_r\vert \right)\, d \vert D_q\mu\vert\\
&+ \int_{B_r} \left(g_r \vert H\vert + \vert f\vert\, \vert D_q g_r\vert\right)\, d \vert D_p\mu\vert\\
&+\int_{B_r\setminus A}  \vert H - D_qf\vert \,\vert D_p g_r\vert \, d\mu\\
&+ \int_{B_r\setminus A}  \vert G - D_pf\vert\,\vert D_q g_r\vert \, d\mu\\
&\leq 
C \big[\vert D_q\mu\vert (B_r) + \vert D_p\mu\vert (B_r)\big]\\
&+\frac{C}{r(1-\rho)}
\big[
\vert D_q\mu\vert (B_r)
+ \vert D_p\mu\vert (B_r)
+ \mu (B_r\setminus A)\big]
\end{split}
 \end{equation*}
 where $C$ is a suitable positive constant independent from $r$ and $\rho$. Consequently, $C$ can be chosen such that we have 
 \begin{equation}
 \label{eq:intBrGgrEST}
 \begin{split}
 \left\vert\int_{B_r}\Gamma g_r\, d\mu\right\vert
&\leq
\frac{C}{r(1-\rho)}
\big[
\vert D_q\mu\vert (B_r)
+ \vert D_p\mu\vert (B_r)
+ \mu (B_r\setminus A)\big],
\end{split}
 \end{equation}
 for all $r,\rho\in (0,1)$. On the other hand
  \begin{equation*}
\begin{split}
\left\vert\int_{B_r}\Gamma g_r\, d\mu\right\vert
&\geq
\left\vert\int_{B_{\rho r}}\Gamma g_r\, d\mu\right\vert 
-
\left\vert\int_{B_r\setminus B_{\rho r}}\Gamma g_r\, d\mu\right\vert
\end{split}
\end{equation*}
 that is
 \begin{equation}
 \label{eq:intBrrGEST}
 \left\vert\int_{B_{\rho r}}\Gamma \, d\mu\right\vert 
 \leq
 \left\vert\int_{B_{r}}\Gamma g_r\, d\mu\right\vert 
 +
 \left\vert\int_{B_r\setminus B_{\rho r}}\Gamma g_r\, d\mu\right\vert 
 \end{equation}
 From \eqref{eq:intBrGgrEST} and \eqref{eq:intBrrGEST} (choosing a larger $C$, if need be), it follows that
 \begin{equation*}
 \begin{split}
\left\vert \frac{1}{\mu (B_{\rho r})}\int_{B_{\rho r}}\Gamma\, d\mu\right\vert
 &\leq 
\frac{C}{r(1-\rho)\mu (B_{\rho r})}
\big[
\vert D_q\mu\vert (B_r)
+ \vert D_p\mu\vert (B_r)
+ \mu (B_r\setminus A)\big]\\
&+\frac{C}{\mu (B_{\rho r})}\, \mu (B_r\setminus B_{\rho r})\\
&=
\frac{C}{1-\rho}\cdot \frac{\mu(B_r)}{\mu (B_{\rho r})}
\left[
\frac{\vert D_q\mu\vert (B_r)}{r\mu(B_r)}
+ \frac{\vert D_p\mu\vert (B_r)}{r\mu(B_r)}
+ \frac{\mu (B_r\setminus A)}{r\mu(B_r)}\right]\\
&+C \left(
\frac{\mu(B_r)}{\mu (B_{\rho r})}-1
\right)
 \end{split}
 \end{equation*}
  for all $r,\rho\in (0,1)$. Hence, by assumptions (iii) and (iv), we obtain
  \begin{equation*}
  \vert D_pH(x) - D_qG(x)\vert
  \leq
  C(\sigma (\rho)-1)
  \end{equation*}
  for every $\rho$ in a left neighborhood of $1$. The conclusion follows from assumption (iii).
 \end{proof}
 
 \begin{remark}
 \label{rem:ThmReducesToSchwarz}
 If $\mu :=\leb^n$, $f\in C^2(\Omega)$, $G:=D_pf$ and $H:=D_qf$ then Theorem \ref{thm:gengenschwarz} reduces trivially to the Schwarz theorem on cross derivatives. However, we cannot claim a new proof of the Schwarz theorem, since the latter was actually used to prove our statement. 
 \end{remark}

  \begin{remark}
  Let us consider a smooth $k$-dimensional surface $S\subset \rn^n$, without boundary or with smooth boundary. Then a hasty attitude might suggest that the distributional derivatives of the Haudorff measure carried by $S$, i.e., $D_i(\hf^k\ree S)$, with $i=1,\ldots,n$, are themselves real Borel measures. Instead, in general this is not the case, and we will show this through the following very simple example. Let $n=2$, $k=1$ and
  $$
  S:=\{(x_1,x_2)\in\rn^2\,\vert\, x_1=x_2\}.
  $$
  Let us set $\mu:=\hf^1\ree S$ for simplicity and observe that
  \begin{equation}
  \label{eq:D_1mExS=diag}
  (D_1\mu)(\varphi)
  =-\int_S D_1\varphi\, d\hf^1 
  =-\sqrt2\int_{\rn}(D_1\varphi) (t,t)\, dt
  \end{equation}
  for all $\varphi\in C^{\infty}_c(\rn)$. Now let $\eta: [0,+\infty)\to [0,1]$ be a decreasing function of class $C^\infty$ such that
  \begin{equation*}
  \label{eq:etainit1then0}
  \eta\vert_{[0,2\pi^2]}\equiv 1,
  \quad
  \eta\vert_{[2\pi^2+1 , +\infty)}\equiv 0
  \end{equation*}
  and define $\varphi_1 , \varphi_2, \ldots \in C_c^{\infty}(\rn^2)$ as follows
  \begin{equation*}
  \varphi_j (x_1 , x_2):=
  \eta (x_1^2+x_2^2)\cos (jx_1) \sin (jx_2).
  \end{equation*}
 From \eqref{eq:D_1mExS=diag} and the equality
  $$
  (D_1\varphi_j)(t,t)
  =2t\eta' (2t^2)\cos (jt) \sin (jt)
  - j\eta (2t^2)\sin^2 (jt)
  $$
  we obtain
  \begin{equation*}
   (D_1\mu)(\varphi_j)
   =-2\sqrt2\, I_j' +
   j\sqrt2\, I_j'',
  \end{equation*}
  with
  \begin{equation*}
  I_j'
  :=\int_{\rn}  t\eta' (2t^2)\cos (jt) \sin (jt) \, dt, 
  \quad
   I_j''
  :=\int_{\rn}  
  \eta (2t^2)\sin^2 (jt)
  \, dt.
  \end{equation*}
  Hence
  \begin{equation}
  \label{eq:StiInfD1m}
  \vert (D_1\mu)(\varphi_j)\vert
  \geq
   j\sqrt2\, \vert I_j''\vert
  -2\sqrt2\,\vert  I_j' \vert 
  =
   j\sqrt2\, I_j''
  -2\sqrt2\,\vert  I_j' \vert 
  \end{equation}
  where
  \begin{equation}
  \label{eq:StiI1}
  \vert I_j'\vert
  \leq \int_{\rn}\vert t\eta' (2t^2)\vert\, dt
  = -2\int_0^{+\infty} t\eta' (2t^2)\, dt
  =-\frac{1}{2}\int_0^{+\infty}D[\eta (2t^2)]\, dt
  =\frac{1}{2}
  \end{equation}
  and
  \begin{equation}
  \label{eq:StiI2}
  I_j''
  \geq\int_{-\pi}^\pi \eta (2t^2)\sin^2 (jt)\, dt 
  =\int_{-\pi}^\pi\sin^2 (jt)\, dt
  =\pi.
  \end{equation}
 From \eqref{eq:StiInfD1m}, \eqref{eq:StiI1} and \eqref{eq:StiI2} we obtain
 \begin{equation}
 \label{eq:FinEstD1m}
 \vert (D_1\mu)(\varphi_j)\vert
  \geq j\pi\sqrt2 -\sqrt2
  \qquad (j=1, 2, \ldots).
  \end{equation}
  Since we have also
  \begin{equation*}
  \label{eq:Estphij}
  \max_{\rn^2}\vert\varphi_j\vert\leq 1,\quad
  \spt\varphi_j\subset B_{2\pi^2+1}(0,0)
   \qquad (j=1, 2, \ldots),
  \end{equation*}
  then the estimate \eqref{eq:FinEstD1m} proves that $D_1\mu$ is not a real Borel measure.
  \end{remark}

 We will now present two simple applications in the context of Lebesgue measure.

 \begin{corollary}
 \label{cor:SchwGenePerhLn}
 Let $h$ be a nonnegative function in $C^1(\rn^n)$. Moreover consider an open set $\Omega\subset\rn^n$, $f, G, H\in C^1(\Omega)$, a couple of integers $p,q$ satisfying $1\leq p< q\leq n$, $x\in\rn^n$ and assume that
 \begin{itemize}
 \item[(i)] $h(x)>0$;
 \item[(ii)] $x\in\Omega \cap A^{h\leb^n , 1}$, where $A$ is the set defined in Theorem \ref{thm:gengenschwarz} (in particular $x$ is in the closure of ${h^{-1} ((0,+\infty))}$);
 \item[(iii)] For $i=p, q$, one has $\int_{B_r(x)}\vert D_ih\vert\,d\leb^n
 = o(r^{n+1})$, as $r\to 0+$ (e.g., $D_ih(y)=o(\vert y-x\vert)$, as $y\to x$).
 \end{itemize}
 Then $D_pH(x)=D_qG(x)$.
 \end{corollary}

\begin{proof}
We will apply Theorem \ref{thm:gengenschwarz} with $\mu:= h\leb^n$. For this purpose, we observe that
 \begin{equation}
 \label{eq:exhLnDim}
 D_i\mu = (D_ih)\leb^n, \text{ hence } \vert D_i\mu\vert = \vert D_ih\vert\leb^n
 \qquad 
 (\text{for all }i=1,\ldots,n)
 \end{equation}
 and (taking into account (i))
 \begin{equation*}
  \label{eq:exhLnsigma}
 \sigma (\rho)=\liminf_{r\to 0+}\frac{\int_{B_r(x)}h\, d\leb^n}{\int_{B_{\rho r}(x)}h\, d\leb^n}=\rho^{-n} \qquad (\text{for all $\rho >0$}).
 \end{equation*}
 Thus assumptions (i), (ii) and (iii) of Theorem \ref{thm:gengenschwarz} are trivially verified. Finally, assumption (iv)  of Theorem \ref{thm:gengenschwarz} is equivalent to (iii) (by (i) and \eqref{eq:exhLnDim}). Therefore Theorem \ref{thm:gengenschwarz} proves the statement. 
\end{proof}

  \begin{remark}
  \label{rem:CorCasoLn}
  If in Corollary \ref{cor:SchwGenePerhLn} we take $h\equiv 1$ then assumptions (i) and (iii) are trivially verified at every $x\in\rn^n$. Recalling also (1) of Remark \ref{rem:PropSuperd}, we conclude that $D_pH=D_qG$ in $\Omega\cap A^{(n+1)}$. In particular, the following property immediately follows:
If $f\in C^1(\Omega)$, $F\in C^1(\Omega , \rn^n)$ and define 
   $A_*:=\{ 
   x\in\Omega\,\vert\, (D_1f (x), \ldots D_nf(x)) = F(x)
   \}$,
  then $DF^t = DF$ in $\Omega \cap A_*^{(n+1)}$.
\end{remark}

\begin{corollary}
\label{cor:teoSchwGenForLnU}
 Let $U\subset\rn^n$ be an open set with boundary of class $C^1$ and let $(\nu_1,\ldots,\nu_n)$ denote the unit outward normal vector field to $\partial U$. Moreover consider $f, G, H\in C^1(\rn^n)$, a couple of integers $p,q$ satisfying $1\leq p< q\leq n$, $x\in\rn^n$ and assume that
 \begin{itemize}
 \item[(i)] $x\in \partial U\cap A^{\leb^n\ree U , 1}$, where $A:=\{ y\in\rn^n \,\vert\, (D_pf (y), D_qf(y)) =(G(y),H(y))\}$;
 \item[(ii)]  For $i=p, q$, one has $\int_{\partial U\cap B_r(x)}\vert \nu_i \vert\,d\hf^{n-1}
 = o(r^{n+1})$, as $r\to 0+$ (e.g., $\nu_i(y)=o(\vert y-x\vert^2)$, as $y\to x$).
 \end{itemize}
  Then $D_pH(x)=D_qG(x)$.
 \end{corollary}
 
 \begin{proof}
 Define $\mu:=\leb^n\ree U$, $\Omega:=\rn^n$ and observe that assumptions (ii) of Theorem \ref{thm:gengenschwarz} is verified by (i), while assumptions (iii) of Theorem \ref{thm:gengenschwarz} follows from
 \begin{equation*}
 \lim_{r\to 0+}\frac{\mu (B_r(x))}{\mu (B_{\rho r}(x))}
 =\lim_{r\to 0+}\frac{\leb^n (U\cap B_r(x))}{\leb^n (U\cap B_{\rho r}(x))}=\rho^{-n}.
 \end{equation*}
 Moreover, by the divergence theorem, we have
 \begin{equation*}
 D_i\mu 
 =
 - \nu_i\,\hf^{n-1}\ree\partial U, \text{ hence } \vert D_i\mu\vert
 =
 \vert\nu_i\vert\,\hf^{n-1}\ree\partial U
 \qquad (i=1,\ldots,n).
 \end{equation*}
  Thus assumption (i) of Theorem \ref{thm:gengenschwarz} is trivially verified, while (ii) yields assumption (iv) of Theorem \ref{thm:gengenschwarz}.
  \end{proof}

\section{Declarations}

\subsection{Ethical Approval} 
This declaration is not applicable.
 
\subsection{Competing interests}
There is no interests of a financial or personal nature.
 
\subsection{Authors' contributions}
This declaration is not applicable (there is only one author).
 
\subsection{Funding}
No funding was received.
 
\subsection{Availability of data and materials}
This declaration is not applicable.

%

 \end{document}